\documentclass[12pt]{amsart}
\usepackage[leqno]{amsmath}

\usepackage{amsmath}
\usepackage{amsthm}
\usepackage{amssymb}
\usepackage{amscd}
\usepackage{amsxtra}     
\usepackage{epsfig}
\usepackage{verbatim}
\usepackage[all]{xypic}

\newtheorem{thm}{Theorem}[section]
\newtheorem{lem}[thm]{Lemma}
\newtheorem{prop}[thm]{Proposition}

\newtheorem{cor}[thm]{Corollary}
\newtheorem{dfn}[thm]{Definition}
\newtheorem{eg}[thm]{Example}
\newtheorem{ques}[thm]{Question}

\theoremstyle{remark}
\newtheorem*{rmk}{Remark}

\numberwithin{equation}{section}

\newcommand{\CC}{\mathbb{C}}

\newcommand{\ZZ}{\mathbb{Z}}
\newcommand{\QQ}{\mathbb{Q}}

\newcommand{\tensor}{\otimes}

 \DeclareMathOperator{\Tor}{Tor}
 \DeclareMathOperator{\Ext}{Ext}
 \DeclareMathOperator{\Hom}{Hom}

 \DeclareMathOperator{\Ass}{Ass}
 \DeclareMathOperator{\Soc}{Soc}
 \DeclareMathOperator{\Supp}{Supp}
 \DeclareMathOperator{\Spec}{Spec}

 \DeclareMathOperator{\Image}{Im}
 \DeclareMathOperator{\Int}{int}
 \DeclareMathOperator{\rad}{rad}
 \DeclareMathOperator{\pd}{pd}
 \DeclareMathOperator{\height}{height}

 \DeclareMathOperator{\depth}{depth}

 \DeclareMathOperator{\syz}{syz}

\newcommand{\Ann}{\textup{Ann}}
 \newcommand{\LC}{\textup{H}}

\begin{document}

\bibliographystyle{plain}

\title{On liftable and weakly liftable modules}
\author{Hailong Dao}
\address{Department of Mathematics, University of Utah,  155 South 1400 East, Salt Lake City,
UT 84112-0090, USA} \email{hdao@math.utah.edu} \maketitle
\begin{abstract}
Let $T$ be a Noetherian ring and $f$ a nonzerodivisor on $T$. We
study concrete necessary and sufficient conditions for a module over
$R=T/(f)$ to be weakly liftable to $T$, in the sense of Auslander,
Ding and Solberg. We focus on cyclic modules and get various
positive and negative results on the lifting and weak lifting
problems. For a module over $T$ we define the loci for certain
properties: liftable, weakly liftable, having finite projective
dimension and study their relationships.
\end{abstract}
\maketitle

\section{Introduction: A brief history of lifting modules}
In  this note, all rings are commutative, Noetherian with
identity, and all modules are finitely generated. Let $T\to R$ be
a ring homomorphism. An $R$-module $M$ is said to {\it lift} (or
{\it litable}) to $T$
 if there is a $T$-module $M'$ if $M = M'
\otimes _T R$ and $\Tor _i^T(M',R)=0$ for all $i>0$. $M$ is said
to  {\it weakly lift} (or {weakly liftable})to $T$  if it is a
direct summand of a liftable module. When $R=T/(f)$ where $f$ is a
nonzerodivisor in $T$, which will be our main focus, then the
$\Tor$ conditions for lifting simply says that $f$ must be a
nonzerodivisor on $M'$. The lifting questions began with:

\begin{ques}(Grothendieck's lifting problem) Let $(T,m,k)$ be a complete
regular local ring and $R = T/(f)$ where $f\in m-m^2$. Does an
$R$-module always lift to $T$ ?
\end{ques}

Note that if $T$ is equicharacteristic, then the answer is
obviously ``yes": in that case $T \cong R[[f]]$, and we can simply
choose $M' = M[[f]]$. The significance of this question was first
publicly realized by Nastold, who observed in \cite{Na} that
Serre' multiplicity conjectures could be solved completely (i.e,
in the case of ramified regular local ring) if we can always lift
in the sense of Grothendieck. Hochster(\cite{Ho1}) gave a negative
answer to Grothendieck's lifting problem (see example \ref{Hoex}).
However, he pointed out that a positive answer to the lifting
problem for prime cyclic modules, and even less would be enough
for Serre's conjectures. Specifically, he posed the following,
which was indeed the starting point for this note:

\begin{ques}\label{Hoques}(Hochster's lifting problem)
Let $(T,m,k)$ be a complete regular local ring and $R = T/(f)$
where $f\in m-m^2$. Let $P\in \Spec(R)$.
\begin{enumerate}
\item When can $M=R/P$ lift ?
\item When there exist an $R$-module
$M$ liftable to $T$ such that $\Supp(M) = \Supp(R/P)$ ?
\end{enumerate}

\end{ques}

Auslander, Ding and Solberg (\cite{ADS})were the first to introduced
and studied systematically the notion of weak lifting. They showed
that in the case $R=T/f$,  weakly lifting an $R$ module $M$ to $T$
is the same as lifting $M$ ``as far as" $T/(f^2)$. If one can repeat
this process to $T$ modulo higher and higher powers of $f$, then one
can lift to $T$ itself, assuming completeness.

Over the years, a number of very interesting results on the lifting
problems have been published. They are almost exclusively
homological in nature. For example, the obstruction to lifting in
Peskine-Szpiro's  thesis can be roughly described as followed : the
fact that $M$ is liftable means that one can lift the whole
projective resolution of $M$ to $T$. This in turn forces certain
module associated to $M$ to have finite projective dimension over
$R$, and that is an obstruction. Using this idea one can construct
modules of finite projective dimension over $R$ (a necessary
condition for liftability when $T$ is regular), but can not lift to
$T$. Jorgensen constructed some very nice examples of such cyclic
modules in \cite{Jo1} (see Example \ref{Joex}). On the positive
side, Buchsbaum and Eisenbud showed that in the case $R=T/(f)$, a
cyclic $R$ module $R/I$ is liftable if  $\pd_R R/I\leq 2$ or if
$\pd_R R/I= 3$ and $I$ is $3$-generated. Jorgensen also produced a
big class of liftable modules, starting from complete intersections
(see \cite{Jo2}).

In this note, we will focus our attention on \textit{concrete}
sufficient and necessary conditions to weak liftability, since
Auslander, Ding and Solberg have made clear that understanding
weak lifting is essential to understanding lifting. Many of our
results are ideal-theoretic, not homological. We have several
motivations for this approach. Firstly, in the context of
Hochster's lifting questions, when $R$ is itself a regular local
ring, if one has to find a negative example, most homological
obstructions would not work ($R$ is ``homologically too nice"). In
any case, to have any hope of answering part (2) of Question
\ref{Hoques} one needs to know `` What annihilates a liftable
module ?". Secondly, for the more general lifting question, it
would be very desirable to tell whether one can  weakly lift a
module just from its presentation. We were able to give some
modest answers to these problems and shed some lights on why they
are non-trivial.

Section 2 reviews basic notations and important results we would
use, including Hochster's characterization of approximately
Gorenstein rings. In Section 3 we study some general necessary and
conditions for weak liftability that involves the annihilator of the
module $M$ (Theorem \ref{theorem3.1}). As applications, we revisit
Hochster's counterexample to Grothendieck's lifting question and
show that it gives a lot more, namely an ideal that is not an
annihilator of any weakly liftable module (see \ref{Hoex}). We also
show that under suitable assumptions, the weakly liftable ideals of
small heights have to be complete intersections (see
\ref{lowheight}).

In Section 4 we focus on weak liftings of cyclic modules. We collect
some simple but useful characterization of weakly liftable cyclic
modules in Lemma \ref{prop4.2}. Many applications follow. We revisit
Jorgensen's example of an unliftable module with finite projective
dimension and give a simple proof in \ref{Joex}, as well as a big
class of such modules in \ref{Joex1}. We also reprove a result
related to modular representation of cyclic groups in \ref{group}. A
negative example to part (1) of Hochster's lifting question above is
given in \ref{Hoex1}. Lastly, we prove very concrete
characterizations of weak liftability for Gorenstein ideals of
dimension $0$ and Cohen-Macaulay, generically Gorenstein ideals of
dimension $1$ in Theorem \ref{thm4.5}.

In Section 5 we formulate a comparative study of liftable, weakly
liftable and finite projective dimension properties. We define a
locus for each property in a quite general way: by fixing a module
over $T$ and asking what hypersurfaces $R$ would make the module
satisfy that property. Our definitions may be viewed as natural
extensions of the notions of ``support sets"  or ``support
varieties" of modules, invented and studied recently by Avramov,
Buchweitz (\cite{AB}) and Jorgensen (\cite{Jo3}). We show in many
cases that weakly liftable and liftable are ``open condition" (see
\ref{prop5.1}, \ref{prop5.2}). This explains in a conceptual way
the existence of many examples of modules with finite projective
dimension but can not lift: they form a Zariski open set in a
certain affine space (see \ref{ex5.3}). Example \ref{ex5.4} and
\ref{ex5.5} show that computing these loci is quite non-trivial,
and in particular the liftable locus may depend on the arithmetic
of the residue field.

Finally, Section 6 contains miscellaneous results and open
questions. We try to emphasize the fact that our knowledge in this
area is still shockingly limited by proposing some simple, yet
intriguing questions.

The author would like to thank Melvin Hochster, whose valuable
insights and advices initiated and inspired most of this work.

\section{Notations and preliminary results}

In  this note, all rings are commutative, Noetherian with
identity, and all modules are finitely generated. Let $R$ be a
ring and $M,N$ be $R$-modules. If $N$ is a submodule of $M$, $N$
is called a \textit{pure} (respectively, \textit{cyclically pure})
if for every $R$-module $E$ (respectively, every cyclic $R$-module
$E$), the induced map $N\tensor E \to M\tensor E$ is injective. If
$M/N$ is of finite presentation, then it is not hard to show that
$N$ is a pure submodule of $M$ is and only if $N$ is a direct
summand of $M$ (see \cite{Ma}, Theorem 7.14).

A more interesting question is when cyclic purity implies purity,
especially when $N=R$. This was answered completely in \cite{Ho2}.
Recall that a local ring $(R,m,k)$ is called \textit{approximately Gorenstein} if for any
integer $N$, there is an ideal $I \subset m^N$ such that $R/I$ is Gorenstein. A Noetherian
ring $R$ is called \textit{approximately Gorenstein} if the localization at any maximal
ideal of $R$ is approximately Gorenstein. Then:
\begin{prop}\label{Ho}(\cite{Ho2}, Proposition 1.4)
Let $R$ be a Noetherian ring. The following are equivalent:\\
(1) $R$ is approximately Gorenstein.\\
(2) For every module extension $R \hookrightarrow M$, cyclic purity
implies purity.
\end{prop}

Hochster's paper also provided very concrete characterizations of
approximately Gorenstein ring. For our purpose, the following
result would be enough:
\begin{thm}(\cite{Ho2}, Theorem 1.7)\label{aG}
Let $R$ be a locally excellent Noetherian ring and suppose that
$R$ satisfies
one of the conditions below:\\
(1) R is generically Gorenstein (i.e., the quotient ring of $R$ is Goresntein).\\
(2) For any prime $P \in \Ass(R)$ and maximal ideal $m \supset P$, $\dim(R/P)_m \geq 2$.\\
Then $R$ is approximately Gorenstein.
 \end{thm}

Let $(R,m,k)$ be a local ring. Let $M,N$ be $R$-modules
such that $l(M \tensor N) < \infty$. One can define the Poincare series
for $M,N$ as :
$$ P^R_{M,N}(t) = \sum_i l(\Tor^R_i(M,N))t^i$$
When $N=k$, we shall simply write $P^R_M(t)$.

The result below is essential for our study of weak lifting.
It is from \cite{ADS} (Proposition 3.2):

\begin{prop} \label{Aus} Consider $R=T/(f)$, where $f$ is a nonzerodivisor on $T$,
which is a Noetherian algebra over a local ring.
The following are equivalent:\\
(1) $M$ is weakly liftable to $T$. \\
(2) $\syz_1^T(M)/f\syz_1^T(M) \cong M\oplus \syz_1^R(M)$, where
$\syz_1^R(M)$ is induced from the free resolution defining
$\syz_1^T(M)$.\\
(3) $M$ is liftable to $R_2 = T/(f^2)$.
\end{prop}

\begin{rmk}
Throughout this paper, when we consider the lifting in the situation $R=T/(f)$, we
will always assume the condition :``T is a Noetherian algebra over a local ring". Since
this covers algebras over fields or DVRs and all local rings, it is not a serious restriction.
\end{rmk}

Finally, we would like to make a definition, mainly for notational
conveniences (see \ref{theorem3.1}).

\begin{dfn}
Let $J,L$ be ideals of a ring $R$. One defines:
\[
    \Int _L(J)  := \{x \in R \ | \exists a_i \in J^i, i=1,..,n
    : \ x^n +a_1x^{n-1} + .. +a_n \in L \}
\]

\end{dfn}

\begin{lem}\label{lem1.4}
It is easy to see that:
\[
\Int _L(J) = \Int _L(J+L) \subseteq \rad (J+L)
\]
\end{lem}

\begin{lem}\label{lem1.5}
If $M$ is a $T$ module and $I = \Ann _T(M)$ then for any ideal $J$ of $T$:
\[
\Ann (M/{JM}) \subseteq \Int _I(J)
\]
\end{lem}

\begin{proof}See \cite{Ma}, Theorem 2.1.
\end{proof}

\section{Some general remarks on weak lifting}

In this section we study several necessary conditions for a module
over $R=T/(f)$ to be weakly liftable to $T$. Our main purpose is
to find concrete obstructions to weak liftability of $M$. Note
that an obstruction to weak lifting is naturally an obstruction to
lifting.

To state the first result, let us recall the change of rings exact
sequence for $\Tor$. Let $R=T/(f)$, where $f$ is a nonzerodivisor on
$T$. Let $M,N$ be $R$-modules. Then we have the long exact sequence
of $\Tor$s :

$$        \begin{array}{ll}
...\to \Tor_{n}^R(M,N) \to \Tor_{n+1}^T(M,N) \to \Tor_{n+1}^R(M,N)\\
\to  \Tor_{n-1}^R(M,N) \to \Tor_{n}^T(M,N) \to \Tor_{n}^R(M,N)\\
\to ...                                                         \\
\to \Tor_{0}^R(M,N) \to \Tor_{1}^T(M,N) \to \Tor_{1}^R(M,N) \to 0

\end{array}
$$\\
In the long exact sequence above, let $\alpha_i$ be the connecting map  $\Tor_{i+2}^R(M,N) \to
\Tor_i^R(M,N)$.

\begin{prop}\label{prop3.1}
Let $T$ be a Noetherian algebra over a local ring. Let $f$ be a
nonzerodivisor in $T$ and $R=T/(f)$. Let $M$ be an $R$-module. The
following are equivalent:
\begin{enumerate}
\item
M is weaky liftable.
\item
The map $\theta : 0 \to M \to \syz_1^TM/f\syz_1^TM$ splits.
\item
For any $R$-module $N$, the map $\alpha_0 : \Tor_2^R(M,N) \to \Tor_0^R(M,N)$ is
$0$.
\item
For any $R$-module $N$  and any integer $i\geq 0$ the map $\alpha_i : \Tor_{i+2}^R(M,N) \to
\Tor_i^R(M,N)$ is $0$.

\end{enumerate}

\end{prop}
\begin{proof}
The equivalence of 1) and 2) is from \cite{ADS}. That 4) implies 3) is obvious.
It remains to show that 2) and 3) are equivalent and 2) implies
4). For that we need to understand how the maps $\alpha_0$ arises. Let:
$$ 0 \to \syz_1^TM \to T^a \to M  $$
be the projective covering of $M$ with respect to $T$. Tensoring
with $R=T/(f)$,since $\Tor_1^T(T,R)=0$ and $\Tor_1^T(M,R)=M$,we get:
$$ 0 \to M \to  \syz_1^TM/f\syz_1^TM \to R^a \to M$$
Breaking down this exact sequence we have:
$$ 0 \to M \to  \syz_1^TM/f\syz_1^TM \to \syz_1^RM \to 0 $$
Tensoring the above exact sequence with $N$ over $R$ gives the connecting map
$\Tor_1^R(\syz_1^RM,N) \to M\tensor_RN  $, which is $\alpha_0$.
From this discussion we can see that 3) is equivalent to the assertion that
the injection $ \theta : M \hookrightarrow \syz_1^TM/f\syz_1^TM$ remains injective when
we tensor with any $R$-module $N$. But this is equivalent to $\theta$ splits (see
\cite{Ma}, theorem 7.14). Also, if $\theta$ splits then all the maps
$\Tor_{i+1}^R(\syz_1^RM,N) \to \Tor_i^R(M,N)$ must also be $0$, which shows
that 2) implies 4).

\end{proof}

The following theorem gives necessary conditions for an ideal to
be the annihilator of a weakly liftable module:
\begin{thm}\label{theorem3.1}
Let $T$ be a Noetherian algebra over a local ring. Let $f$ be a
nonzerodivisor in $T$ and $R=T/(f)$. Let $M$ be an $R$-module and
$I=\Ann_T(M)$.
If $M$ is weakly liftable to $T$ then:\\
1) $(I^2:f) \subseteq I$ \\
2) $(JI:f) \subseteq \Int _I(J)$ for all ideals $J$ of $T$\\
3) $(JI:f) \subseteq \rad (I+J)$ for all ideals $J$ of $T$
\end{thm}

We begin with some lemmas. Let us try to understand concretely
what weak liftability imposes on the annihilator of a module. Let
$M$ be an $R$-module and we pick a free covering of $M$ as a
$T$-module:
\[
      0 \to W \to G \to M \to 0
\]

Here $G=T^n$. Let $I = \Ann _T(M)$.

By the above Proposition, the map $\theta$:
\[
\xymatrix{0 \ar[r] &  G/W \ar[r]^{h} & W/fW }
\]
which takes $x+W$ to $fx+fW$ splits.

\begin{lem}\label{lem3.2}
Let $T,R,M,G,W$ as above. If $M$ is weakly liftable to $T$ then for
any ideal $J \subseteq T$:
\[
(JW:f) \subseteq (JG+W)
\]
\end{lem}

There are two proofs of this lemma. The first is very elementary.
The second enables us to apply Hochster's results to strengthen
the conclusions in the cyclic case(see next section).
\begin{proof} (proof 1)
By Proposition \ref{prop3.1} \ $G/W \cong fG/fW$ is a direct summand of $W/fW$. So there is
a submodule $B$ of $W$ such that:\\
 1) $fW \subseteq B \subseteq W$ \\
 2) $B+fW=W$ \\
 3) $B \cap fG \subseteq fW$  \\
Now suppose $v \in (JW:f)$. So there are $j_i$'s in $J$ and $w_i$'s in $W$
such that : $fv = \sum{j_iw_i}$ . But from 2) each $w_i = fg_i+b_i$ with
$g_i \in G$ and $b_i \in B$. So we have:
\[
fv = \sum{j_i(fg_i+b_i)}
\]
Rearranging:
\[
f(v- \sum{j_ig_i}) = \sum{j_ib_i}
\]
Since LHS is in $fG$ and RHS is in $B$, from condition 3) we get
$v-\sum{j_ig_i} \in W$ or $v \in (JG+W)$.
\end{proof}

\begin{proof} (proof 2)
We use the simple fact that for $T$-modules $P \subseteq  Q$ such that $P$
is a direct summand of $Q$ , then for any ideals $J$ of $T$ , $P/JP$ injects
into $Q/JQ$  (in other words, $P$ is a cyclically pure submodule of $Q$).\\
Applying that to $G/W$ and $W/fW$ we have $G/(W+JG)$ injects into $W/(fW+JW)$
 ( with the map induced from $h$), which is equivalent to :
\[
  (fW+JW) : f \subseteq (W+JG)
\]
which can be easily seen to be equivalent to :
\[
   (JW:f) \subseteq (W+JG)
\]
\end{proof}

\begin{lem} Let $T,R,M,G,W,I$ be as above. Then for any
ideal $J$ in T :
\[
(JI:f) \subseteq \Ann _T(G/(JW:f))
\]
\end{lem}

\begin{proof}
Let $v \in (JI:f)$ . So $vf \in JI$ . Hence $vfG \subseteq JIG$ . But $I$ kills
$G/W$ , so $IG \subseteq W$. It implies that $vfG \subseteq JW$ $\Longrightarrow$
$vG \subseteq (JW:f)$ $\Longrightarrow$ $v \in \Ann _T(G/(JW:f))$.
\end{proof}

Now we can prove Theorem \ref{theorem3.1}:
\begin{proof}(of \ref{theorem3.1})
By the previous Lemmas we have : \\
$(JI:f) \subseteq (\Ann _T(G/(JW:f)) \subseteq
\Ann _T(G/(JG+W))
= \Ann _T((G/W)/(J(G/W))) = \Ann _T(M/JM)$ \\ The last term is $I$
if $J=I$ , and it is contained in $\Int _I(J)$ otherwise (by \ref{lem1.4}).
Finally, by \ref{lem1.5} we have $\Int _I(J) \subseteq \rad (I+J)$, as required.
\end{proof}
As an application we will revisit Hochster's counterexample to
Grothendieck lifting question (see \cite{Ho1}).
\begin{eg}\label{Hoex}
Let $T=\ZZ_{(2)}[[x,y,z,a,b,c]]$. Let
$f=2$ and $R=T/(f)$.  Let $I=(2,x^2,y^2,z^2,a^2,b^2,c^2,xa+yb+zc)$
and  $g=xayb+ybzc+zcxa$. Because of the relation :
\[
   2g=(xa+yb+zc)^2 - x^2a^2+ y^2b^2+ z^2c^2
\]
It follows that $g \in (I^2:f)$. But is is not hard to show $g \notin I$.
By \ref{theorem3.1}, not only $T/I$ is not liftable to $T$, as Hochster showed, but
$I$ can not be the annihilator of any $R$-module which is  weakly liftable to $T$.
\end{eg}

Let $R=T/(f_1,f_2,..,f_c)$ where the $f_r$'s form a $T$-sequence. Then
the definition of liftability and weak liftability is unchanged. Note that the
condition $\Tor _i^T(M',R)=0$ for all $i>0$ is equivalent to the $f_r$'s form
a regular $M'$-sequence. It is probably worth mentioning:

\begin{cor}
Let $T,R,M$ and $f_r$'s as above. Let $I=\Ann_T(M)$. Suppose M is weakly
liftable to $T$. Then for each $1 \leq r \leq c$ : \\
(1) $(I^2:f_r) \subseteq I$ \\
(2) $(JI:f_r) \subseteq \Int_I(J)$ for all ideals $J$ of $T$ \\
(3) $(JI:f_r) \subseteq \rad (I+J)$ for all ideals $J$ of $T$ \\

\end{cor}

\begin{proof}
We only need to prove for $f_1$. Suppose $M$ is a direct summand of
$M_1$, which lifts to $M_2$, a $T$-module. Then viewed as a
$T/(f_1)$ module, $M_1$ lifts to $M_2/(f_2,..,f_c)$. So $M$, as
$T/(f_1)$-module, is weakly liftable. Now we only need to apply
Theorem \ref{theorem3.1}.
\end{proof}

Next, we present another simple corollary of \ref{theorem3.1}:
\begin{cor}\label{poincare}
Let $(T,m.k)$ be a local ring and $R=T/(f)$ where $f$ is a nonzerodivisor in $T$.
Suppose $M,N$ are $R$-modules such that $M\tensor N$ is of finite length and $M$ is
weakly liftable to $T$. Then $P_{M,N}^T(t)=(t+1)P_{M,N}^R(t)$. If $T$ is regular,
$M$ is weakly liftable to $T$ and $\dim M<\dim R$, then $(t+1)^2 \mid P_M^T(t)$.
\end{cor}

\begin{proof}
By Theorem \ref{theorem3.1}, the change of rings long exact
sequence for $\Tor$ would break down into short exact sequences:
$$ 0 \to \Tor_{i}^R(M,N) \to \Tor_{i+1}^T(M,N) \to \Tor_{i+1}^R(M,N) \to 0 $$
for all $i\geq 0$. The first statement is immediate. As for the second, first note that
$\pd_R M<\infty$. Since $\dim M<\dim R$, $P_M^R(-1) = \chi^R(M,k) = 0$. So $(t+1)\mid P_M^R(t)$,
this fact and the first statement finish the proof.

\end{proof}

As an application, we will show that weakly liftable Cohen-Macaulay or Gorenstein ideals of small heights
often are complete intersections:

\begin{cor}\label{lowheight}
Let $(T,m.k)$ be a regular local ring and $R=T/(f)$ where $f$ is a nonzerodivisor in $T$.
Let $I$ be an ideal in $R$ such that $R/I$ is weakly liftable to $T$. If $\height(I)=1$ and
$R/I$ is Cohen-Macaulay then $I$ is principal. If $\height(I)=2$ and
$R/I$ is Gorenstein then $I$ is generated by two elements.
\end{cor}

\begin{proof}
Let $J$ be the preimage of $I$ in $T$. By Corollary \ref{poincare} we have  $(t+1)^2 \mid P_{T/J}^T(t)$.
In the first case $P_{T/J}^T(t)$ has to be equal to $(t+1)^2$ (because $\pd_T T/J=2$). In the second
case $P_{T/J}^T(t)$ has to be equal to $(t+1)^3$(because $\pd_T T/J=3$ and the last Betti number is $1$
since $T/J$ is Gorenstein). In both cases we must conclude that $J$ is a complete intersection,
and so is $I$.
\end{proof}

\begin{eg}
Let $T =k[[x_1,...,x_n]]$, $f=x_1$ and $R=k[[x_2,...,x_n]]$. Then any $R$-module is liftable
to $T$ and the above corollary says that in $R$, a height $1$ Cohen-Macaulay ideal has to
be principal and a height $2$ Gorenstein ideal has to be $2$-generated. So there is little hope
to strengthen the result.
\end{eg}

\section{Weakly liftable cyclic modules}

In the case of cyclic modules, the statements of the previous
section can be simplified or strengthened. Let us recall the basic
setup. Let $T$ be a Noetherian  algebra over a local ring and $f$ be
a nonzerodivisor in $T$. Let $R=T/(f)$ and $I$ be an ideal in $T$
which contains $f$. We will focus on finding conditions for $T/I$ to
be weakly liftable (as an $R$-module) to $T$.

\begin{lem}\label{prop4.2}
Let $T,f,R,I$ be as above. Fix $\bold v = (f,f_1,..,f_n)$ a set
of generators for $I$. The following are equivalent:\\
(1) $M=T/I$ is weakly liftable to $T$.\\
(2) The $T$-linear map $h: T/I \to I/fI$ which takes $1+I$ to $f+fI$
splits.\\
(3) The $T$-linear map $g: T/I \to I/I^2$ which takes $1+I$ to $f+I^2$
splits.\\
(4) For any presentation of $I$:
$$\xymatrix{T^{m} \ar[r]^X & T^{n+1} \ar[r]^{\bold v} & I \ar[r] & 0}$$
Let $\bold r,\bold r_1,...,\bold r_n $ be the rows of $X$. There exist
$x_1,...,x_n \in T$ such that :
$$ \bold r - x_1\bold r_1+...+x_n\bold r_n \in IT^m $$

And they imply the following equivalent conditions :\\
(5) $(IJ:f) \subseteq (J+I)$ for any ideal $J$. \\
(6) $(IJ:f) \subseteq J $ for any ideal $J \supseteq I$. \\
(7) (If $T$ is local) $(IJ:f) \subseteq J $ for any irreducible ideal $J$.\\
If in addition, $T/I$ is approximately Gorenstein, then all
the conditions (1) to (6) (and (7) in the local case) are equivalent.
\end{lem}

\begin{rmk}
The last assertion (when $T/I$ is approximately Gorenstein)  was
first suggested in \cite{Ho1}, page 462.
\end{rmk}

\begin{proof}
The equivalence of (1) and (2) is a restatement of \ref{theorem3.1}.
If (2) holds, then $I/fI = T/I\oplus N$ for some $T$-module $N$.
Tensoring with $T/I$ we get : $I/I^2 = T/I \oplus N/IN$, which
gives (3). Now assume (3) which says the map $g$ splits.  But $g$ is a
composition of
$$\xymatrix{T/I \ar[r]^h & I/fI \ar[r] & I/I^2 }$$
 so $h$ also splits.

For the equivalence of (3) and (4), let $Z = \Image (X)$ be the first syzygy
of $I$. Tensoring the exact sequence :
$$ 0 \to Z \to T^{n+1} \to I \to 0$$
with $T/I$ we get:
$$ 0 \to (Z\cap IT^{n+1})/IZ \to Z/IZ \to (T/I)^{n+1} \to I/I^2 \to 0$$
which shows that $Z/(Z\cap IT^{n+1})$ is a first syzygy of $I/I^2$
(as a module over $T/I$). So there is no new relations, and
$I/I^2$ admits the following presentation:
$$\xymatrix{\overline{T}^{m} \ar[r]^{\overline X} & \overline{T}^{n+1} \ar[r]^{\overline{\bold v}} & I/I^2 \ar[r] & 0}$$

Here $\bar{}$ denotes mod $I$. Then (3) means exactly that there exist
$x_1,...,x_n \in T$ such that :
$$ \overline{\bold r} = \overline{x_1} \overline{\bold r_1}+...+\overline{x_n} \overline{\bold r_n} $$
Next, (1) implies (5) is a restatement of Lemma \ref{lem3.2}. The
equivalence of (5) and (6) is trivial.
The only thing to check now is equivalence of (6) and (7). Clearly (6)
implies (7). Suppose (6) fails and we have an ideal $J$ such that
$(IJ:f) \nsubseteq J $. Pick $x \notin J$ such that $xf \in IJ$.
Choose a maximal ideal $J_1$ containing $J$ such that $x\ \notin
J_1$. Then $J_1$ is irreducible, and (7) fails as well.

Finally, suppose that in addition $T/I$ is approximately Gorenstein.
Condition (4) says that the map $g$, viewed as a $T/I$-module extension,
is cyclically pure. Then Proposition \ref{Ho} implies that $T/I$ is a pure submodule
of $I/I^2$ via $g$, so (3) holds. That finishes our proof.
\end{proof}

\begin{eg}
We give an example to show that if $T/I$ is not approximately
Gorenstein, the last assertion of Lemma \ref{prop4.2} would fail
even in simplest cases. Let $T=\QQ[[x,y]]$, $m=(x,y)$,$I=m^2$ and
$f=x^2+y^2$. Clearly $Im:f\subset m$ and $I^2:f\subset I$. Let $J$
be any ideal lying strictly between $I$ and $m$. Then
$J=m^2+(ux+vy)$, with $u,v\in \QQ$. We want to show that $IJ:f
\subset J$. Pick $g\in m$ such that $fg \in IJ = m^4+(ux+vy)m^2$.
Let $g'$ be the linear part of $g$, then clearly $g'f\in
(ux+vy)m^2$. Since $f$ is irreducible in $T$,$g' \in (ux+vy)$,
thus $g\in J$. So condition (6) of Lemma \ref{prop4.2} is
satisfied. However $T/I =T/m^2$ is not weakly liftable to $T$. One
can see it by using Theorem \ref{Joex1} or simply observing that
$\pd_{T/(f)} T/m^2 = \infty$.
\end{eg}

It is now quite easy to show that one of the main examples in a paper by Jorgensen
(example 3.3 in \cite{Jo1} ) gives a cyclic module of finite projective dimension but is
unliftable:

\begin{eg}\label{Joex}
Let $k$ be a field, $T=k[[x_1,x_2.x_3,x_4]]$, $f=x_1x_2-x_3^2$, $R=T/(f)$,
$I = (f,i_1,i_2,i_3,i_4)$, where: \\
$i_1 = -x_2x_3+x_2x_4 ,\ i_2 = x_1x_3+x_2x_3,\
i_3 = -x_2^2-x_3x_4 ,\ i_4 = x_1^2-x_2^2+x_3^2-x_4^2$\\
Finally, let $J =(x_1,x_3,x_4,x_2^2) \supset I$. It can be shown using
Macaulay that $\pd _R T/I = 3$. But $-x_3b_1+x_4b_2+x_1b_3 = x_2f$, so
$x_2 \in (JI:f)$. Obviously $x_2 \notin J$, so $T/I$ is not even weakly
liftable.
\end{eg}

The above example suggests the following:
\begin{thm}\label{Joex1}
Let $T=\oplus_{n\geq 0}T_n$ be a graded ring with $T_0=k$ is a field. Let $I$ be a $T$-ideal
generated by homogeneous elements of degree $a$. Let $f \in I$ be a homogeneous nonzerodivisor of degree $a$
such that $(f) \subsetneq I$.
Assume that $I$ admits a free presentation:
$$\xymatrix{F \ar[r]^{X} & G \ar[r]^Y & I \ar[r] & 0}$$
such that all the entries of the matrix $X$ has degree $b<a$. Then
$T/I$ as a module over $R=T/(f)$ is not weakly liftable to $T$.
\end{thm}

\begin{proof}
As $f$ must be a $k$-linear combination of the generators of $I$, we may
as well assume that $Y = (f,f_1,...,f_n)$. Then let $\bold r,\bold r_1,...,\bold r_n $
be the rows of $X$. By part (4) of \ref{prop4.2} there exist
$x_1,...,x_n \in T$ such that :
$$ \bold r - x_1\bold r_1+...+x_n\bold r_n \in IT^m $$
Counting degree, there must be $y_1,...,y_n \in k$ such that
$$ \bold r = y_1\bold r_1+...+y_n\bold r_n $$
But this means that $(f)$ is a direct summand of $I$ as $T$-modules.
This is impossible unless $(f)=I$, so we are done.

\end{proof}
As another application, we would prove the following, which is
relevant to the theory of modular representation of cyclic groups
(see \cite{The}). We give a brief explanation. Let $D$ be a
discrete valuation ring whose maximal ideal is generated by a
prime number $p$. Let $C_p$ be the cyclic group of order $p$. Let
$A =D/p^2$ and $k=D/pD$. One wishes to study the $AC_p \cong
A[X]/(X^p-1)$-modules. Let $M$ be such a module. Then $M/pM$ is a
$kC_p \cong k[X]/(X^p-1)\cong k[X]/(X-1)^p$ module. The
decomposable modules over $kC_p$ must be of the form
$S_i=k[X]/(X-1)^i$. So $M/pM$ is a direct sum of $S_i$'s. The
interesting questions is which $i$ may occur ? Clearly this
corresponds to when is $S_i$  liftable to $AC_p$, or equivalently,
weakly liftable to $DC_p$ (by \ref{Aus}). In view of this, the
following corollary is a special case of Theorem 5.5 in \cite{The}
:

\begin{cor}\label{group}
Let $(D,m,K)$ be a discrete valuation ring whose maximal ideal is generated by
a prime number $p$. Let $T = D[X]/(X^p-1)$, $R = T/(p) \cong K[X]/(X^p-1)\cong K[X]/(X-1)^p$.
Let $S_i = K[X]/(X-1)^i$ ($1\leq i\leq p$) be $R$-modules. Then $S_i$ is weakly liftable to
$T$ is and only if $i \in \{1, p-1,p\}$.
\end{cor}

\begin{proof}
Clearly $S_p=R$ lifts and $S_1$ lifts (take $S=T/(X-1)$, then $S$ is a lift
of $S_1$. We assume $1<i<p$. Note that $S_i = T/(p,(X-1)^i)$. Over $T$, the
ideal $I = (p,(X-1)^i$ has a presentation:
$$\xymatrix{T^{2} \ar[r]^X & T^{2} \ar[r]^{\bold v} & I \ar[r] & 0}$$
Here $\bold v =(p,(X-1)^i)$ and $X$ has 2 rows: $\bold r = ((X-1)^p,
g(X))$ where $g(X) = \frac{(X^p-1)-(X-1)^p}{p}$ and $\bold r_1 =
(-p, (X-1)^{p-i})$. By Theorem \ref{prop4.2} (equivalence of (1) and
(4)), $T/I$ is weakly liftable if and only if $g(X)$ is a multiple
of $(X-1)^{p-i}$ (mod $I$). Rewriting:
$$g(X) = \frac{((X-1+1)^p-1)-(X-1)^p}{p} = \sum_{j=1}^{p-1}\frac{ \binom{p}{j}}{p} (X-1)^j$$
One can see that it happens if and only if $p-i=1$.
\end{proof}

Next we gives an example in which $R=T/(f)$ is a ramified regular
local ring of dimension 11 and a prime cyclic module of $R$ that
is not weakly liftable. This shows that there is a negative
example to part (1) of Question \ref{Hoques}.
\begin{eg}\label{Hoex1}
Let $T=V[[x,y,z,a,b,c,u,v,w,t]]$, in which $(V,2V)$ is a  DVR. Let
$f=2$ and $R=T/(f)$ and let $\bar{}$ denote mod $f$. Abusing
notation, we don't use $\bar{}$ for the indeterminates. Let
$I=(2,tu-x^2,tv-y^2,tw-z^2,xa+yb+zc)$. Since $t$ is not nilpotent
modulo $I$, we can pick a minimal prime $\overline P$ over
$\overline I$ which doesn't contain $\overline t$. It is easy to see
that actually, $\overline P=\overline I:t^\infty = \overline I: t$.
Using Macaulay 2, we can actually calculate
$P=(I,ua^2+vb^2+wc^2,uayz+vbzx+wcxy)$. For our purpose, we only need
to see that $\overline P \subset (u,v,w,x^2,y^2,z^2,xa+yb+zc)$. Now,
let $P$ be the preimage of $\overline P$ in $T$,
$J=(P,t,a^2,b^2,c^2)$ and $g=xayb+ ybzc+zcxa$. Because of the
relation :
\[
   2g=(xa+yb+zc)^2+(tu-x^2)a^2+(tv-y^2)b^2+(tw-z^2)c^2-t(ua^2+vb^2+wc^2)
\]

It follows that $g \in PJ$. It suffices to show $g \notin J$. We
can do so modulo $2,u,v,w,t$. Then because of remark above, it is
enough to show $g \notin (x^2,y^2,z^2,a^2,b^2,c^2,xa+yb+zc)$. But
this is true by \ref{Hoex}. By \ref{prop5.2} we can replace
$f$ by $2+f'$ with $f' \in mP$ to get an example where $R$ is an
honest ramified regular local ring.
\end{eg}
\begin{rmk}
Similar examples surely exist for all characteristics.
\end{rmk}

Lemma \ref{prop4.2} still leaves much to be desired when one wants
to show some module to be weakly liftable, since checking cyclic
purity involves infinitely many ideals $J$. To really take
advantage of the conditions, we need a few lemmas:

\begin{lem}\label{lem4.3}
Let $(T,m,k)$ be a local ring and $I \subseteq J_1 \subseteq J_2$
be ideals in $T$. Assume that $T/J_1$ is $0$-dimensional and
Gorenstein (in other words, $J_1$ is irreducible). Then $IJ_2:f
\subseteq J_2$ if $IJ_1:f \subseteq J_1$.
\end{lem}

\begin{proof}
Suppose the assertion is not true. Then we can find $x$ such that
$fx \in IJ_2$ but $x \notin J_2$. Since $T/J_1$ is Gorenstein and
$0$-dimensional, $\Hom(-,R/J_1)$ is a self-dualizing functor. As
$J_2 \subsetneq J_2+(x)$ we must have $\Hom(T/(J_2+(x)),T/J_1)
\cong J_1:(J_2+(x)) \subsetneq J_1:J_2 \cong \Hom(T/J_2, T/J_1)$.
So we can pick $y \in J_1:J_2$ but $y \notin J_1:(J_2+(x))$. Then
$fxy \in IJ_2y \subset IJ_1$. By assumption this forces $xy \in
J_1$ which implies $y(J_2+(x)) \subset J_1$, contradicting our
choice of $y$.
\end{proof}

\begin{lem}\label{lem4.4}
Let $(T,m,k)$ be a local ring and $I \subseteq J$ be ideals in
$T$. Assume that $T/J$ is $0$-dimensional and Gorenstein. Let
$u\in T$ represent the generator of the socle of $T/J$. Then $IJ:f
\subset J$ if and only if $fu \notin IJ$.
\end{lem}

\begin{proof}
One direction is clear, so assume $\ uf \notin IJ$ and let $J_1 =
IJ:f$. If $J_1 \supsetneq J$ then we let $n$ to be the smallest
integer such that $m^nJ_1 \subseteq J$. By assumption $n \geq 1$
and $m^{n-1}J_1 \nsubseteq J$. Let $s \in m^{n-1}J_1$ but $s
\notin J$. Then $ms \subseteq J$, so $as-u \in J$ for some unit
$a$. But $s$ is clearly in $J_1$ (here we need $n\geq 1$), so $fs
\in IJ$. Then $fu = fs + f(s-u) \in IJ$, a contradiction.
\end{proof}

\begin{thm}\label{thm4.5}

Let $R=T/(f)$ where $(T,m,k)$ is a local ring and $f$ is a
nonzerodivisor in $T$. Let $T/I$ be an $R$-module (so $f\in I$).\\
(1) Suppose that $T/I$ is $0$-dimensional and Gorenstein. Let $u\in
T$ represent the generator of the socle of $T/I$. Then $T/I$
is weakly liftable if and only if $\ uf \notin I^2$. \\
(2) Suppose that $T/I$ is  1-dimensional, Cohen-Macaulay and  generically
Gorenstein. Let $J \subset T$ represent the canonical ideal of
$T/I$. Let $u\in T$ represent the generator of the socle of $T/J$.
Then $T/I$ is weakly liftable if and only if $\ uf \notin IJ +
I^{(2)}$.

\end{thm}

\begin{proof}(1) By Lemma \ref{prop4.2} and Lemma \ref{lem4.3}.\\
(2) Let $ S = T/I$. Then since $S$ is generically Gorenstein, its
canonical module $\omega_S$ is isomorphic to an ideal of height
$1$. Let $\overline J$ be that ideal in $S$ (here $J$ is an ideal
in $T$ and $\bar{}$ denotes modulo $I$. We claim that $S/\overline
J$ is $0$-dimensional and Gorenstein. Since $\overline J$ is
height $1$, the first assertion is trivial. Now apply $\Hom(k,-)$
to the short exact sequence :
$$ 0 \to \overline J \to S \to S/\overline J $$
and observe that $\Hom(k,S) = 0$ since $\depth S = 1$ we get:
$$ 0 \to \Hom(k,S/\overline J) \to \Ext^1_R(k,\overline J)$$
Since $\overline J \cong \omega_S $ we can use local duality to
get $\Ext^1_R(k,\overline J) \cong \Ext^1_R(k,\omega_S) \cong
\LC^0_m(k)^{\vee} \cong k$. So $\Hom(k,S/\overline J)$ injects
into $k$ and since it is not zero, it has to be $k$. So
$S/\overline J$ is Gorenstein. Let $\overline x$ be a
nonzerodivisor in $S$. Then $x\overline J \cong J \cong \omega_S$
so $x\overline J$ must also be an irreducible ideal. Note that
$xu$ represent the generator of $\Soc(S)$. By Lemma \ref{prop4.2}
and Lemma \ref{lem4.3} we only need to check that $xuf \notin
I(I+xJ)$ for any $x$ such that $\overline x$ is a nonzerodivisor
in $S$ . This is equivalent to $uf \notin IJ + (I^2:x)$ for all
such $x$, or $uf \notin IJ + I^{(2)}$ as desired.

\end{proof}

\section{The (non) liftable and weakly liftable loci}
This section is a comparative study of liftable, weakly liftable
and finite projective dimension properties. Throughout the section
we will assume that $(T,m,k)$ is a  local ring, and $M$ is a
$T$-module. Let $I \subset \Ann_T(M)$ be an ideal in $T$ and fix a
minimal system of generators $(f_1,...,f_n)$ for $I$ . Then there
is a map $\alpha : I \to k^n \cong I/mI$ induced by
$(f_1,...,f_n)$. For a property $\mathcal{P}$ we define the
$\mathcal{P}$-locus of $M$ in $I$ as :
$$\mathcal{L}_{\mathcal{P}}(I,M) := \{ f\in I \vert M \ \text{satisfies} \ \mathcal{P} \ \text{as a module
over} \ T/(f) \}  $$ and the geometric $\mathcal{P}$-locus of $M$
in $I$ as :

$$\mathcal{V}_{\mathcal{P}}(I,M) := \alpha(\mathcal{L}_{\mathcal{P}}(M))$$

If $I= \Ann_T(M)$ we shall simply write
$\mathcal{L}_{\mathcal{P}}(M)$ and $\mathcal{V}_{\mathcal{P}}(M)$.
For $\mathcal{P} = \{\text{not liftable}\}$ (resp. not weakly
liftable, not finite projective dimension) we will write
$\mathcal{L}_{nl}$ (resp. $\mathcal{L}_{nwl}$,
$\mathcal{L}_{npd}$) (by convention $0$ is in all of these sets)
and $\mathcal{V}_{nl}$ (resp. $\mathcal{V}_{nwl}$,
$\mathcal{V}_{npd}$). It is more convenient to work with the
negative properties, as they turns out to be ``closed" conditions.

\begin{rmk}
When $(f_1,...,f_n)$ form a regular sequence on $T$, then
$\mathcal{V}_{\mathcal{P}}(I,M)$ agrees with the ``support
variety" of $M$ as defined in \cite{AB}. When $I=\Ann(M)$,
$\mathcal{V}_{\mathcal{P}}(M)$ agrees with the ``support set" of
$M$ defined in \cite{Jo3}.
\end{rmk}
We first observe that :
\begin{prop}\label{compare}
Suppose $T$ is a regular local ring and $M$ is a $T$-module. Let
$I=\Ann_T(M)$. Then:
$$ I \supset \mathcal{L}_{nl}(M) \supset \mathcal{L}_{nwl}(M) \supset \mathcal{L}_{npd}(M) \supset mI $$
and
$$ k^n \supset \mathcal{V}_{nl}(M) \supset \mathcal{V}_{nwl}(M) \supset \mathcal{V}_{npd}(M) $$
\end{prop}
\begin{proof}
The only thing needs to be proved is $\mathcal{L}_{npd}(M) \supset
mI $. Let's assume $f \in mI$ and $R=T/(f)$. By a result of
Shamash (\cite{Sha}), in this situation:
$$P^T_M(t) = (1-t^2)P^R_M(t)$$
which clearly shows that the $P^R_M(t)$ can not be finite series
(otherwise $P^T_M(t)$ would have negative terms!).
\end{proof}

\begin{prop}\label{prop5.1}
$\mathcal{L}_{nl}(T/I)$ is an ideal.
\end{prop}
\begin{proof}
First, let $f\in \mathcal{L}_{nl}(T/I)$ and $a\in T$. We want to show $af \in \mathcal{L}_{nl}(T/I)$.
Assume it is not true, so there exists a $T$-ideal $J$ such that $af$ is a nonzerodivisor on $\overline T = T/J$
and $J+(af)=I$. The first condition shows that $f$ is also a nonzerodivisor on $\overline T$, and the second
shows that $\overline f\overline T \subset \overline a \overline f\overline T$. By Nakayama's Lemma, $\overline a $ is an unit in $\overline T$, so
$\overline T$ is also a lift of $T/I$ with respect to $f$.

Secondly, let $f,g \in \mathcal{L}_{nl}(T/I)$. Similarly, suppose $f+g \notin  \mathcal{L}_{nl}(T/I)$, we
seek a contradiction. Again, there there exists a $T$-ideal $J$ such that $f+g$ is a nonzerodivisor on $\overline T = T/J$
and $J+(f+g)=I$. Since $f,g \in I$ we must have, in $\overline T$, $\overline f = (\overline f + \overline g)e_1$ and
$\overline g = (\overline f + \overline g)e_2$. Adding the two equations and using that $f+g$ is a nonzerodivisor on $\overline T$, we
get $e_1+e_2=1$ in $\overline T$. This forces $e_1$ or $e_2$ to be a unit in $\overline T$, but then $\overline T$ must be a lift
of $T/I$ with respect to either $f$ or $g$.
\end{proof}

\begin{prop}\label{prop5.2}
If $T/I$ is approximately Gorenstein, then $\mathcal{L}_{nwl}(T/I)$ is an ideal.
\end{prop}
\begin{proof}
We first construct a sequence $\{L_i\}$ of irreducible ideals in $T/I$ such that
$L_{i+1} \subsetneq L_i \ \forall i$ and  $\{L_i\}$ in $T/I$ are cofinal with the powers
of the maximal ideal in $\overline T = T/I$. Just pick $L_1$ as any irreducible ideal in $\overline T$.
Then there is a power of $\overline m$, $\overline m^l \subset L_1$. By assumption we
can pick an irreducible ideal $L_2 \subset \overline m^l$, and so on.
Let $J_i$ be the preimage of $L_i$ in $T$. By \ref{prop4.2} and \ref{lem4.3}
$f\in \mathcal{L}_{nwl}(T/I)$ if and only if $IJ_i:f\nsubseteq J_i$ for some $i$ (since any
irreducible ideal would contain some $J_i$). Let $I_i := \{f\in I \vert IJ_i:f\nsubseteq J_i\}$.
By \ref{lem4.4} $I_i = (IJ_i:s_i)\cap I$, here $s_i$ represent the socle element of $J_i$.
So each $I_i$ is an ideal in $T$. But \ref{lem4.3} and the fact that $J_{i+1}\subseteq J_i$
shows that $I_i \subseteq I_{i+1}$. Hence the sequence of ideals $\{I_i\}$ must stabilize, and since
$\mathcal{L}_{nwl}(T/I) = \cup_1^{\infty} I_i$  we are done.
\end{proof}

\begin{eg}\label{ex5.3}
Proposition (\ref{prop5.1}) implies that $\mathcal{V}_{nl}(T/I)$
is an affine space. So as long as $\mathcal{V}_{npd}(T/I)$ is not
a linear algebraic set, then there should be quite a few example
of finite projective dimension, unliftable cyclic modules: they
form the non-empty Zariski open set
$\mathcal{V}_{nl}(T/I)\setminus \mathcal{V}_{npd}(T/I) $ in
$\mathcal{V}_{nl}(T/I)$. Such nonlinear $\mathcal{V}_{npd}(T/I)$
are known to be quite common, see the examples at the end of
\cite{Jo3}.
\end{eg}

\begin{eg}\label{ex5.4}
Theorem \ref{thm4.5} gives explicit formula for
$\mathcal{L}_{nwl}(T/I)$ in some cases. Specifically, using the
notations of Theorem \ref{thm4.5} we have $\mathcal{L}_{nwl}(T/I)=
I^2:u$ when $T/I$ is  Gorenstein of dimension $0$ and
$\mathcal{L}_{nwl}(T/I)= (IJ + I^{(2)}):u$ if $T/I$ is
Cohen-Macaulay, generically Gorenstein of dimension $1$.
\end{eg}

\begin{eg}\label{ex5.5}
Let $T=k[[X,Y,Z]]/(X^2+Y^2+Z^2)$, here $k$ is a field. Let $x,y,z$
be the images of $X,Y,Z$ respectively and let $m=(x,y,z)$. We
claim that $\mathcal{L}_{nl}(T/m) = m^2$ if $k=\CC$ and
$\mathcal{L}_{nl}(T/m) = m$ if $k=\QQ$.

First, let $k=\CC$. Choose any element $f = ax+by+cz$ with $a,b,c\in
\CC$. We have to show $f\notin \mathcal{L}_{nl}(T/m)$, in other
words, $T/m$ is liftable to $T$ as a $T/(f)$-module. Let
$I_1=(x,y+iz), I_2=(y,z+ix),I_3=(z,x+iy)$. Note that they are prime
ideals of height 1 in $T$. We claim that one of these ideals
together with $f$ will generate $m$. Let $V_i =\alpha(I_i)$ (so for
example $V_1$ is generated by the vectors $(1,0,0)$ and $(0,1,i)$).
Then the planes $V_1,V_2,V_3$ intersect at only the origin in
$\CC^3$ so one of them, say $V_1$, can not contain the vector
$(a,b,c)$. This shows that $(I_1,f)=m$. But $f$ is clearly a
nonzerodivisor on $T/I_1$, and so $T/m$ is liftable.

Next, assume $k=\QQ$. It suffices to show that $x\in
\mathcal{L}_{nl}(T/m)$, as then $y,z \in \mathcal{L}_{nl}(T/m)$ by
symmetry and hence $m=(x,y,z) \subseteq \mathcal{L}_{nl}(T/m)$ by
Proposition \ref{prop5.1}. Suppose $T/I$ is a lift of $T/m$ as a
module over $T/(x)$. Then $I+(x)=m$. So there are $a,b\in T$ such
that $y-ax, z-bx \in I$. But since $x^2+y^2+z^2=0$ this forces
$x^2(1+a^2+b^2) \in I$. Since $k=\QQ$, $(1+a^2+b^2)$ must be a
unit, hence $x^2\in I$. But then $x$ can not be a nonzerodivisor
on $T/I$.

Finally, observe that $\mathcal{L}_{nwl}(T/m) = m^2$ in both
cases. Indeed, by the previous example, since the socle element of
$T/m$ is $1$, we have $\mathcal{L}_{nwl}(T/m) = m^2:1 = m^2$.

\end{eg}

\section{Miscellaneous results and open questions}

In this section we first collect some  observations relevant to
Grothendieck's lifting question. We begin by noting that in this
case, condition (3) of Theorem \ref{theorem3.1} (the weakest
obstruction) is of no value:

\begin{prop}
Suppose $(T,m,k)$ is a regular local ring and $f\in m-m^2$. Then
for any ideal $I,J$ of $T$: $(IJ:f) \subseteq \rad (I+J)$.
\end{prop}

\begin{proof}
Let $v \in (IJ:f)$ Let $P$ be any prime containing $I+J$. We want
to show that $v \in P$. Localize at $P$ we see that $v \in
(I_PJ_P:f) \subseteq ((P_P)^2:f) $. But $f$ is also a regular
element of $T_P$, so that implies $v \in P$.
\end{proof}

The following result explains why in example \ref{Hoex1}, one
needs $f$ to involve only nonlinear monomials of the
indeterminates:

\begin{prop}
Suppose $T=V[[x_1,..,x_n]]$, where $V=(V,pV,k)$ is a DVR. Let
$\bar{}$ denotes modulo $p$. Let $f \in T$ be such that $\overline
f \in \overline m - \overline m ^2$. Suppose $P$ is a prime ideal
containing $p,f$. Then $T/P$ is weakly  liftable to $T$ as an
$T/(f)$ module.
\end{prop}

\begin{proof}
Suppose $T/P$ is not weakly liftable.  Since $T/P$ is a complete
local domain, it is approximately Gorenstein by (\ref{aG}). So
condition (6) of (\ref{prop4.2}) there is an ideal $J \supset P$
and $v \notin J$ such that $fv \in JP$. Working mod $p$ we have a
counter example in the ring $\overline T = k[[x_1,..,x_n]]$ whose
maximal ideal is $\overline m$ and $\overline f = g \in \overline
m - {\overline m}^2$. But in this case, $\overline T/ \overline P$
is liftable to $T$(in fact, any module is), a contradiction.
\end{proof}

It is natural to ask whether we could obtain some  obstructions
for the class of modules with finite projective dimension over $R$
similar to Theorem (\ref{theorem3.1}). Obviously, we expect such
obstructions to be weaker, since weak liftablility implies finite
projective dimension. In deed, in the example of Hochster([Ho1])
$(I^2:f) \subseteq I$ fails, but $\pd_R T/I$ is still finite
because $R$ is regular.  Surprisingly, the obstruction (3)  still
works:

\begin{prop}
Let $(T,m,k)$ be a regular local ring,$f\in m$ ,$M$ and $R$-module
and $I=\Ann_T(M)$. Suppose $\pd_RM < \infty$. Then for any ideal
$J$ of $T$, $(JI:f) \subseteq \rad(I+J)$.
\end{prop}

\begin{proof}
We only need to prove $(PI:f) \subseteq P$ for a prime $P
\supseteq I$. Suppose this fails for some $P$. Localize at $P$ we
get $f \in P_PI_P$ in the local ring $T_P$. But $\pd_{R_P}M_P$ is
still finite and $R_P=T_P/(f)$, contradicting Proposition
\ref{compare}.
\end{proof}

\begin{eg}
In example (\ref{Joex}) we have $\pd_R T/I < \infty$ and
 $(JI:f) \ni x_2 \notin J$. Note that, however, $x_2 \in \rad(J)$.
\end{eg}

Finally, we would like to pose some questions. Keeping up with the
theme of this note, they are concrete and hopefully realistic:

1) In the situation of Grothendieck's (or Hochster's) lifting
question, is there an example of weakly liftable but not liftable
module ? The same question can be asked even when $(T,m,k)$ is a
regular local ring and $f$ any nonzero element (so $f$ could be in
$m^2$). Example \ref{ex5.5} shows there are plenty of examples
when $T$ is not regular.

2) Can one get necessary conditions for liftability stronger than
those in Theorem \ref{theorem3.1} ? This is vital to have any hope
of answering completely Hochster's question (\ref{Hoques}).

3) Are $\mathcal{L}_{nl}(M)$ and  $\mathcal{L}_{nwl}(M)$ ideals ?
Are there explicit formulas (or algorithms) to compute them ? Or
at least, the dimensions of $\mathcal{V}_{nl}(M)$ and
$\mathcal{V}_{nwl}(M)$ (assuming they are vector spaces)?

4) Under what conditions $\mathcal{V}_{nl}(M)$ (or
$\mathcal{V}_{nwl}(M)$) would be the linear closure of
$\mathcal{V}_{npd}(M)$ ?

\end{document}